\documentclass[conference]{IEEEtran}
\usepackage[bookmarks=false]{hyperref}
\usepackage{amsmath,amssymb,amsfonts}
\usepackage{algorithmic}
\usepackage{graphicx}
\usepackage{textcomp}
\usepackage{hyperref}
\usepackage{multirow}
\usepackage{subcaption}
\usepackage[usenames,dvipsnames]{xcolor}
\usepackage[english]{babel}

\def\BibTeX{{\rm B\kern-.05em{\sc i\kern-.025em b}\kern-.08em
    T\kern-.1667em\lower.7ex\hbox{E}\kern-.125emX}}

\newdimen\slantmathcorr
\def\oversl#1{
\setbox0=\hbox{$#1$}
\slantmathcorr=\wd0
\hskip 0.2\slantmathcorr \overline{\hbox to 0.8\wd0{%
\vphantom{\hbox{$#1$}}}}
\hskip-\wd0\hbox{$#1$}
}

\def\undersl#1{
\setbox0=\hbox{$#1$}
\slantmathcorr=\wd0
\underline{\hbox to 0.8\wd0{%
\vphantom{\hbox{$#1$}}}}
\hskip-0.8\wd0\hbox{$#1$}
}    

\newcommand{\norm}[1]{\left\lVert#1\right\rVert}

\graphicspath{{figures/},.}

\title{Collective Effects and Performance of Algorithmic Electric Vehicle Charging Strategies}

\author{\IEEEauthorblockN{Miroslav~Gardlo\IEEEauthorrefmark{1}, \v{L}ubo\v{s}~Buzna\IEEEauthorrefmark{1}, Rui Carvalho\IEEEauthorrefmark{2}, Richard Gibbens\IEEEauthorrefmark{3} and Frank Kelly\IEEEauthorrefmark{4}}
\IEEEauthorblockA{\IEEEauthorrefmark{1}University of \v{Z}ilina, Univerzitn\'{a} 8215/1, \v{Z}ilina, Slovakia, Emails: gardlo.miroslav@gmail.com, lubos.buzna@fri.uniza.sk} 
\IEEEauthorblockA{\IEEEauthorrefmark{2}Department of Engineering, Durham University,\\ Lower Mountjoy, South Road, Durham, DH1~3LE, UK, Email: rui.carvalho@durham.ac.uk}
\IEEEauthorblockA{\IEEEauthorrefmark{3}Department of Computer Science and Technology, University of Cambridge,\\  William Gates Building, 15 JJ Thomson Avenue, Cambridge, CB3~0FD,~UK, Email: richard.gibbens@cl.cam.ac.uk}
\IEEEauthorblockA{\IEEEauthorrefmark{4}Statistical Laboratory, Centre for Mathematical Sciences, University of Cambridge,\\ Wilberforce Road, Cambridge, CB3~0WB,~UK, Email: f.p.kelly@statslab.cam.ac.uk}}

\begin{document}
\maketitle

\begin{abstract}
We combine the power flow model with the proportionally fair optimization criterion to study the control of congestion within a distribution electric grid network. The form of the mathematical optimization problem is a convex second order cone that can be solved by modern non-linear interior point methods and constitutes the core of a dynamic simulation of electric vehicles (EV) joining and leaving the charging network. The preferences of EV drivers, represented by simple algorithmic strategies, are  conveyed to the optimizing component by real-time adjustments to user-specific weighting parameters that are then directly incorporated into the objective function. The algorithmic strategies utilize a small number of parameters that characterize the user's budgets, expectations on the availability of vehicles and the charging process. We investigate the collective behaviour emerging from individual strategies and evaluate their performance by means of computer simulation.
\end{abstract}

\begin{IEEEkeywords}
EV charging strategies, optimal power flow model, proportional fairness
\end{IEEEkeywords}

\section{Introduction}
Electricity congestion problems may become a barrier to large-scale adoption of EVs. A typical EV requires around 1 kWh for every 3-4 miles of driving, potentially making EVs significant energy consumers in the future~\cite{Mukherjee_2015}. Intelligent management of EV charging can reduce the need for expensive new generation, transmission, and distribution facilities by shifting and controlling the load demand caused by EVs. Since the state of the electric network evolves in time, owners of EVs would have to supervise charging of their vehicle in real time and make decisions on a scale of minutes, which is definitely impossible. A solution is to develop algorithmic strategies (software agents) that take responsibility for decisions that are associated with the charging of EVs. Replacement of human decision making by software agents is already notable in areas as stock markets, personal assistants or autonomous vehicles. 

There has been a large research activity in the last few years in the area of EV charge scheduling with the major efforts summarized in survey papers. Studies analysing the impacts of EVs on distribution networks considering aspects such as driving patterns, charging characteristics, charge timing, and EV penetration are summarized in~\cite{Green_2010}. A survey of studies analysing the impacts of EVs in a broader sense (grid impacts, economic and environmental impacts) was given in~\cite{Richardson_2013}. A review paper of EV charging scheduling approaches~\cite{Mukherjee_2015} distinguishes  two large classes: unidirectional and bidirectional charging. Each class is further classified based on whether the scheduling is centralized or distributed and whether the mobility aspects, renewable energy sources and  ancillary services (e.g. frequency control, generation control) are considered. Furthermore, the paper overviews the type of objective function and the optimization techniques that were used. Other review papers address EV charging from specific perspectives such as grid management~\cite{Mahmud_2018}, the interaction of EVs and smart grids~\cite{Mwasilu_2014}, the optimization techniques used~\cite{Tan_2016} and EV modelling~\cite{Daina_2017}.

The mathematical concept of proportional fairness as a network congestion mechanism in electricity networks has been considered already by several authors. In~\cite{Carvalho_2015}, proportional fairness has been shown as superior to max-flow, while evaluating the critical arrival rate of vehicles leading to the congestion on the electricity network. An EV charging approach based on proportional fairness while considering linear and additive flow capacity constraints, and using a decentralized solving approach based on dual prices was proposed in~\cite{Ardakanian_2013}. Work~\cite{Fan_2012} builds on~\cite{Kelly98,Gibbens_1999} and proposes a decentralized scheme to allocate proportionally fair charging rates to EVs. The network is represented as a price signal that informs the charging scheme about the congestion level. Willingness to pay and price determine the charging rate following the framework proposed in~\cite{Kelly98}. This paper discusses the main principles, but the interaction of users is not investigated and the choice of the willingness to pay parameter is limited to constant values or random numbers. In this paper, we further extend this body of literature by combining the concept of proportional fairness with an optimization model that takes into account the network in a more realistic way and by extending a mathematical description of algorithmic charging strategies~\cite{Knapp_2016} (see Section~\ref{sec_mathematical formulation}). In Section~\ref{sec_results}, we present results of selected numerical experiments and evaluate them. Section~\ref{sec_conclusions} concludes the paper by summarizing the main findings. 

\section{Problem formulation}
\label{sec_mathematical formulation}

\subsection{Concept of the proposed approach}
\label{subsec_concept}

EV charging is an ideal venue for demand response as the charging load is more elastic compared to other residential loads (e.g. cooking, heating etc.). Often, users park their vehicles in the evening  when returning from work and need them to be ready again in the morning. We investigate an approach that exploits the elasticity of EV loads to distribute available network capacity to EVs while the line voltage levels remain within a predefined range. We control the charging of EVs in real-time. In our model, the interests of a user are represented by an algorithmic entity or software agent that acts on their behalf. EV drivers only select the strategy and set a few parameter values. Depending on the allocated budget, state of charge or departure time, software agent $l$ determines the value of the willingness to pay parameter $w_l(t)$ and sends it to the electric power distribution network. As only a single value $w_l(t)$ is revealed to the network the majority of a user's sensitive data is protected.  The network collects $w_l(t)$ values from all users and combines the Optimal Power Flow model with the Proportional Fairness criterion to determine the values of power $P_l(t)$ that are allocated to vehicles.
  
\subsection{Optimization problem}
\label{subsec_opt_model}

A mathematical optimization problem forms the core of the simulation model. We denote the set of vehicles that are willing to pay a positive amount for electric energy at time $t$ as $\mathcal{N}(t)$. The solution of an optimization model is used to determine the vector of electric powers $P(t) = (P_l(t): l \in \mathcal{N}(t))$ that are allocated to electric vehicles, while maintaining the parameters of the power grid network within operational limits. The feasible set of power allocations is denoted as $\mathcal{P}(t)$. We solve the optimization problem 
\begin{IEEEeqnarray}{rCl}
\label{eq:model1_obj_a}
\mbox{$\underset{ P(t)}{\text{maximise}}$} \quad & \sum\limits_{ l \in \mathcal{N}(t)} w_l(t) \log(P_l(t)) \\
\label{eq:model1_con_b}
\text{subject to } \quad & P(t) \in \mathcal{P}(t),
\end{IEEEeqnarray}
where $w_l(t)$ are non-negative weights. Problem (\ref{eq:model1_obj_a})-(\ref{eq:model1_con_b}) can be interpreted as a network protocol that distributes network capacity to vehicles. For $w_l(t)=1$, we recover the proportional fairness allocation~\cite{Kelly98}. However, by including values $w_l(t)$ as weights, the shares can be arbitrarily changed according to a weighted proportionally fair allocation. Values $P_l(t)$ are non-negative and when $w_l(t) = 0$, then the term with $P_l(t)$ is excluded from~(\ref{eq:model1_obj_a}) and set to value of $0$. 

The mathematical description of the set of feasible solutions $\mathcal{P}(t)$, is derived from the optimal power flow problem~\cite{TaylorBook15}. For this purpose, we model the distribution network by a directed tree graph $\mathcal{G}(\mathcal{V}, \mathcal{E})$, where $\mathcal{V}$ is a set of nodes (called buses) and $\mathcal{E}$ is set of edges (called branches). By the symbol $\mathcal{V}^{+}$, we denote the set of all nodes excluding the root node (feeder). Voltage drops are significant in distribution networks, and to a large extent determine the network capacity, which leads us to using the model of power flow specific to distribution networks~\cite{Kersting01}. 

To derive the mathematical description of the set of feasible solutions $\mathcal{P}(t)$, we applied the procedure described in~\cite{Low14b,Low14c}. In the first step, power flows are expressed in terms of complex nodal voltages $v_i(t)$, where $i$ is a network node. To avoid quadratic expressions, nodal voltages for network edges connecting node $i$ with the node $j$ are
substituted by complex variables $V_{ij}(t) = v_i(t)v_j(t)^{*}$ using the symbol $^{*}$, to denote the complex conjugate.
The resulting problem is non-convex, but if the network is radial and if electric powers extracted or injected at nodes and electric powers flowing through power lines are not constrained, we can apply SOCP relaxation that is exact~\cite{TaylorBook15} and leads to the optimization model~(\ref{eq:model4_obj_a})-(\ref{eq:model4_con_f}).
\begin{figure*}[!t]
\begin{align}
\label{eq:model4_obj_a}
\mbox{$\underset{V(t), P(t)}{\text{maximise}}$} \quad & \sum\limits_{l \in \mathcal{N}(t)} w_l(t) \log(P_l(t)) \\
\text{subject to } \nonumber \\ 
\label{eq:model4_con_b}
& \sum\limits_{j: e_{ij} \in \mathcal{E}} (g_{ij}(V_{ii}(t) - Re\{V_{ij}(t)\}) + b_{ij}Im\{V_{ij}(t)\}) = \sum_{l \in C_i(t)} P_l(t) & \qquad i \in \mathcal{V}^{+} \\
\label{eq:model4_con_c}
& \sum\limits_{j: e_{ij} \in \mathcal{E}} (b_{ij}(V_{ii}(t) - Re\{V_{ij}(t)\}) - g_{ij}Im\{V_{ij}(t)\}) = 0 & \qquad i \in \mathcal{V}^{+} \\
\label{eq:model4_con_d}
& \undersl{v_i}^2 \leq V_{ii}(t)
\leq  \oversl{v_i}^2 & \qquad i \in \mathcal{V} \\
\label{eq:model4_con_f}
& \norm{\begin{pmatrix}
2Re\{V_{ij}(t)\} \\
2Im\{V_{ij}(t)\} \\
V_{ii}(t) - V_{jj}(t)\end{pmatrix}}_2 \leq V_{ii}(t) + V_{jj}(t) & \qquad e_{ij} \in \mathcal{E}.
\end{align}
\end{figure*}
The quantities $g_{ij}$, and $b_{ij}$ denote the conductance and susceptance of power lines, respectively. If $X$ is a complex quantity, then its real (active) part is denoted as $Re(X)$ and imaginary (reactive) part is $Im(X)$. Constraints (\ref{eq:model4_con_b}) and (\ref{eq:model4_con_c}) ensure the conservation of active and reactive power at individual network nodes, respectively. Constraints (\ref{eq:model4_con_d}) keep the square of the voltage magnitude within the lower bound $\undersl{v_i}^2$ and upper bound $\oversl{v_i}^2$. Constraints~(\ref{eq:model4_con_c}) are a by-product of the substitution $V_{ij}(t) = v_i(t)v_j(t)^{*}$, maintaining the specific relation between the variables  $V_{ii}(t)$, $V_{jj}(t)$ and $V_{ij}(t)$.

\subsection{Charging strategies}
\label{subsec_user_strategies}

Charging strategies are inspired by the digital file transmission strategies that were originally proposed for use in communication networks~\cite{Gibbens_1999}. The charging strategy is an algorithm updating the willingness to pay parameter, $w_l(t)$, of a vehicle $l$ in order to influence the charging process.  All the presented strategies assume a limited budget $W_{max}^l$, such that the time integral of $w_l(t)$ cannot exceed this value. Moreover, to facilitate the comparison between strategies, it is assumed that vehicles have a limited time $T_{max}^l$ to stay connected to the charger. Consequently, a connected vehicle $l$ is disconnected from the network in three cases: its battery is fully charged, the budget $W_{max}^l$ allocated to the charging operation was completely spent, or the maximum connected time $T_{max}^l$ elapsed.
 
\subsubsection{Uniform spending in time (UT strategy)}
As a reference strategy, we introduce this simple strategy, where the parameter $w_l^{UT}(t)$ is set to a constant value $\frac{W_{max}^{l}}{T_{max}^{l}}$ throughout the charging process. Thus, in this case the whole budget is spent uniformly over time during the time period $T_{max}^{l}$ that the vehicle $l$ is charging.

\subsubsection{Uniform charging in time (UC strategy)}
Vehicles represented by this strategy aim to gain a constant amount of energy per unit of time, while using the entire time interval allocated for charging. Hence, vehicle $l$ updates its willingness to pay parameter using the formula 
\begin{equation}
\resizebox{0.42\textwidth}{!}{$
w_{l}^{UC}(t) =  \max \left\{ 0,  w_{l}^{UC}(t - \Delta t) - \kappa\Delta t\left(B_l(t) - \frac{t - T_{arr}^l}{T_{max}^l}B_{max}^{l}\right) \right\}$},
\end{equation}
where $T_{arr}^l$ is the time when the vehicle $l$ was plugged into the network. The UC strategy increases its willingness to pay parameter $w_l^{UC}(t)$, when the state of battery $B_l(t)$ is lower than the target value $\frac{t - T_{arr}^l}{T_{max}^{l}} B_{max}^{l}$  and decreases $w_l^{UC}(t)$ otherwise. Thus, a user expects that the battery should be charged linearly in time, while reaching full capacity $B_{max}^{l}$ at the time of departure. A limited budget may lead to spending the budget before the time interval $T_{max}^{l}$ has elapsed. Here and throughout $\kappa>0$ is a scaling parameter.

\subsubsection{Affordable spending (AF strategy)}
This strategy is inspired by the file-transfer strategy introduced in~\cite{Gibbens_1999}. On average, the algorithm tries to pay a price $W_{max}^{l}/B_{max}^{l}$ per unit of energy. The user $l$ updates its willingness to pay by comparing the current value $w_l^{AF}(t)$ paid per unit of received power $P_l(t)$ with the price that this user can afford to pay per unit of energy $\frac{W_{rem}^l(t)}{  B_{max}^l -  B_l(t)}$, where $W_{rem}^l(t)$ is the remaining budget that vehicle $l$ can still use for charging at time $t$. The willingness to pay parameter $w_l^{AF}(t)$ is updated as follows:
\begin{equation}
\resizebox{0.42\textwidth}{!}{$w_l^{AF}(t) = \max\left\{ w_l^{AF}(t - \Delta t) + \kappa \Delta t\left(\frac{W_{rem}^l(t)}{B_{max}^{l} - B_l(t)} - \frac{w_l^{AF}(t - \Delta t)}{P_l(t)}\right), w_{min} \right\}$}.
\end{equation}
Parameter $w_{min} > 0$ allows the algorithm to restart the charging process after any time periods when it appears too expensive to charge.
A possible drawback of this strategy is that it does not consider the time limit $T_{max}^{l}$, and thus it takes no action as the departure time of a vehicle is approaching.

\subsubsection{Combination of affordable spending with uniform spending in time (AFT strategy)}
Strategy AFT builds on the AF strategy and makes it sensitive to the time limit $T_{max}^l$. The motivation behind this strategy is to maximize the gained energy when the time limit is approaching even if not saving any budget. Strategy AFT either pays the affordable amount per unit of energy or towards the end of the time limit $T_{max}^l$, it pays a fraction of what is affordable to pay per unit of time. A point at which the strategy starts depleting the budget is controlled by the pre-factor $\alpha(t) \in (-\infty, 1]$, i.e.
\begin{equation}
w_i^{AFT}(t) = \max\left\{w_l^{AF}(t), \alpha(t)\frac{W_{rem}^{l}(t)}{T_{max}^{l} - t}\right\}.
\label{eq_AFT_1}
\end{equation}
We chose
\begin{equation}
\alpha(t) = \frac{t - T_{arr}^l}{T_{max}^{l}(1-d)} - \frac{d}{1-d}, 
\label{eq_AFT_2}
\end{equation}
where the parameter $d \in [0,1)$ defines the point (expressed as the fraction of $T_{max}^l$), when the linearly growing function $\alpha(t)$ is equal to zero. 

The optimization model (\ref{eq:model4_obj_a})-(\ref{eq:model4_con_f}) was solved with CVXOPT solver (\url{http://cvxopt.org}) and the simulation model was implemented in Python programming language and is available for download from~\url{https://bit.ly/2qO9CkX}.

\section{Numerical experiments}
\label{sec_results}
Our results are organized in two sections.  In section~\ref{sec_microscopic_scenarios}, we analyse the behaviour of charging strategies at the microscopic level, in a pre-designed stylized situation, to explore mutual interactions of vehicles and collective phenomena. In section~\ref{macroscopic_performance}, the average performance of the charging strategies at the macroscopic level is evaluated.

\subsection{Exploration of collective properties of algorithmic strategies at microscopic level}
\label{sec_microscopic_scenarios}

\begin{figure}[!ht]
\centering
\includegraphics[scale=0.65]{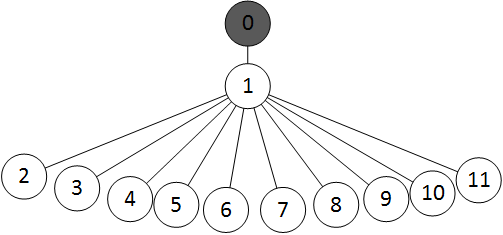}
\caption{Network used in the scenario to explore the collective properties of the algorithmic strategies at microscopic level.}
\label{collective_effects_scenarios}
\end{figure}

Preliminary experiments  have shown that the network topology has an important effect on power allocations. Firstly, we eliminate the network effects from the behaviour of strategies by using a network with the structure displayed in Figure~\ref{collective_effects_scenarios}. Vehicles connect only to nodes $2-11$ and each line $e_{ij} \in \mathcal{E}$ has resistance $r_{ij} = 0.1$ and reactance $x_{ij} = 0.6$ (these are average values estimated from a real-world network~\cite{Gan_2015}). Thus, all vehicles compete for the capacity of the edge that is connecting nodes $0$ and $1$. Furthermore, in all the experiments presented in this paper constraints~\eqref{eq:model4_con_d}, take  the form
\begin{equation}
    \label{eq:model_4_con_g}
    V_{nominal}^2 (1-\alpha)^2 \leq  V_{ii} \leq V_{nominal}^2 (1+\alpha)^2 \quad i \in \mathcal{V},
\end{equation}
where $V_{nominal} = 1.0, \alpha = 0.1$ and we set $d = 0.75$ and $\kappa = 0.001$.

We considered a population of vehicles with heterogeneous budgets and arrival of "aggresive" vehicle. In this scenario exactly one vehicle is connected to nodes $l~=~2,~\dots,~10$, with $B_{max}^l~=~20$, $T_{max}^l~=~300$, $B(0)~=~0$ and  $W_{max}^l~=~1.3^{l-2}500$. The node $11$ is initially empty. At time $t~=~100$ an aggressive vehicle following UT strategy is connected to node $11$ with budget $W_{max}^{11}~=~10~000~000$, $B_{max}^{11}~=~20$ and $T_{max}^{11}~=~100$. The purpose of this experiment is to illustrate the behaviour of strategies when vehicles posses different budgets and the network experiences the arrival of aggressive vehicle possessing a huge budget that needs to be charged fast. Since the ”aggressive vehicle” will consume the majority of the resources of the electric network, this experiment also models a situation, where the electric network moves suddenly into a congested state, in which the total demand of individual vehicles cannot be fulfilled.
\begin{figure}[ht!]
\centering
\includegraphics[scale=0.15]{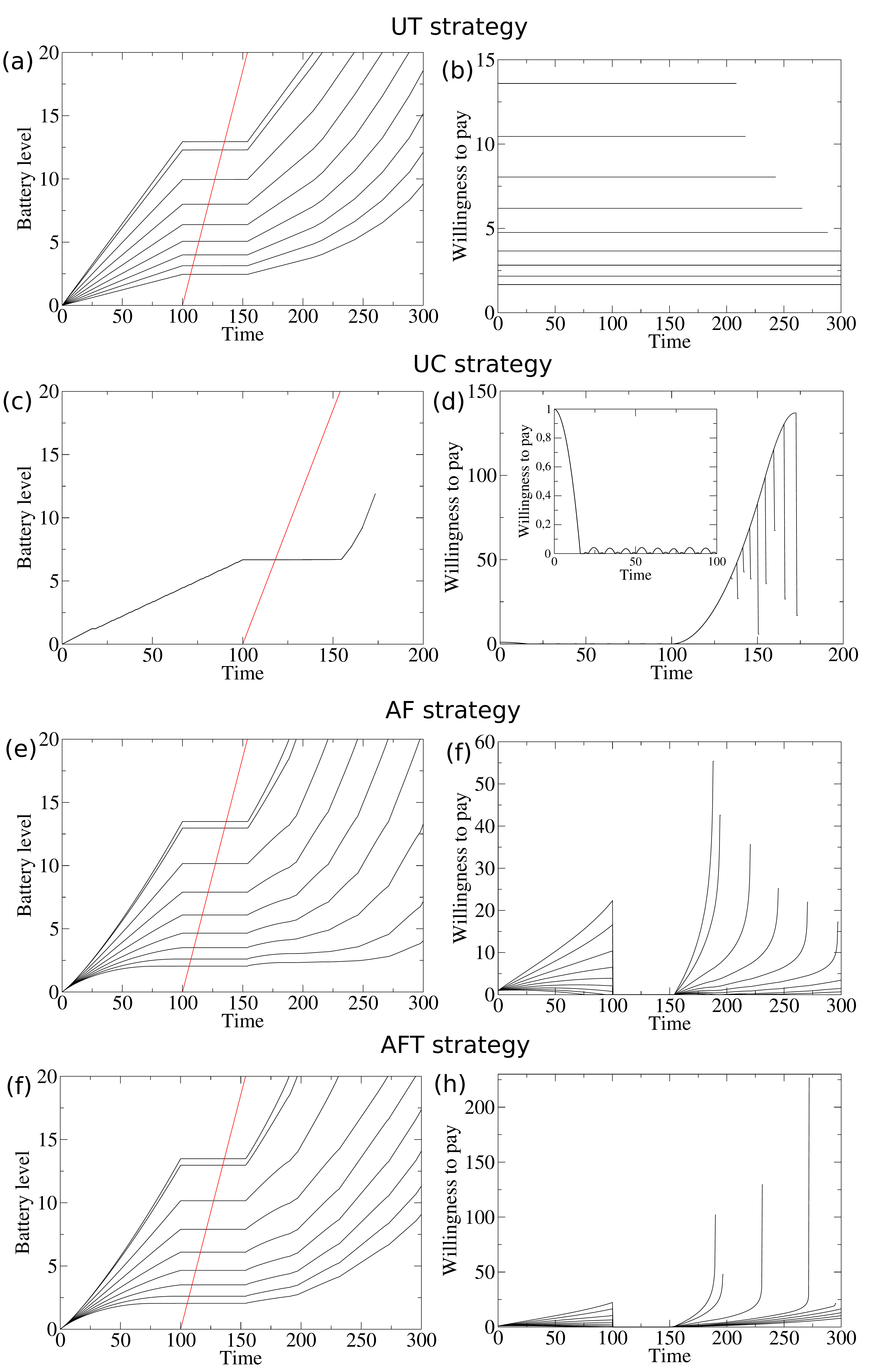}
\caption{Evolution of battery state and willingness to pay parameters for a scenario with heterogeneous budgets and the arrival of aggressive vehicle. We plot zoom of the first 100 time units of the willingness to pay evolution in the inset of panel (d). An aggressive user that possesses huge budget and follows UT strategy arrived to the network at time $t = 100$. The evolution of the state of battery for aggressive user is shown in red colour.}
\label{fig:_aggressive}
\end{figure}
Vehicles using the UT strategy, are charged in the non-increasing order of the available budget, thus the  vehicle charging with the highest rate is the vehicle with the maximum budget. Since with this strategy the users do not update their willingness to pay parameter during the charging process, the only factor impacting the speed of charging process is the number of vehicles in the network. Figure~\ref{fig:_aggressive}a shows that charging process of vehicles with lower budgets speed up as vehicles with higher budgets leave the network. After the arrival of aggressive vehicle, other users following the UT strategy  continue to receive electric energy at a much slower rate. The time limit $T_{max}^l~=~300.0$  was insufficient to charge fully the four vehicles that started with the smallest budgets.

In the Figure~\ref{fig:_aggressive}c, it is not possible to distinguish individual vehicles. While updating willingness to pay in time $t$, the UC strategy only considers the current battery state $B(t)$ and the remaining time in the network $T_{max}^{l} - t$. Since these quantities are independent of the remaining budget, the charging process of all connected vehicles is the same. Initially, vehicles do not need to compete over the capacity and consequently vehicles consume only a very small part of the budget. The inset of Figure~\ref{fig:_aggressive}d zooms the region of the first 100 time units to better see the evolution of the willingness to pay.  As soon as the aggressive vehicle arrived, other vehicles start to compete for electric power by increasing the willingness to pay parameter, because their charging processes did not fulfil demanded constant intake of electric energy. This competition resulted in spending the entire budget of vehicles before they could reach full battery capacity (see Figure~\ref{fig:_aggressive}d).

Vehicles following the AF strategy are charged in the same order as for the UT strategy. Even though the order is the same, there is one distinct behavioural difference between the two strategies, which can be best seen on the vehicles that are charged with the lowest speed. The sorting effect is more pronounced for AF strategy as smaller number of vehicles is being charged simultaneously. Since vehicles with higher budget can afford to pay more, vehicles with lower budget set their willingness to pay to $w_{min}$ and wait for the price to decrease. After the arrival of the aggressive vehicle, all charging vehicles detected that the price per unit of energy is too high, and they lowered their willingness to pay to the lowest possible value $w_{min}$. This can be seen in the Figure~\ref{fig:_aggressive}f. When the aggressive vehicle leaves the network, the price per unit of energy decreases and remaining vehicles increase their willingness to pay parameter again. This behaviour is quite reasonable, as it is essentially impossible to compete against the aggressive user. In this case, three vehicles did not manage to charge fully.

The behaviour of AFT strategy is in this experiment qualitatively similar to the AF strategy. The sorting effect and waiting for lower price per energy unit can be also observed. Compared to the AF strategy only some subtle differences appear that can be seen in Figure~\ref{fig:_aggressive}h, where the increase in the willingness to pay parameter towards the end of the simulation is higher than in the Figure \ref{fig:_aggressive}f. Close examination of Figures~\ref{fig:_aggressive}e and~\ref{fig:_aggressive}f reveals that strategy AFT charges fully one vehicles less. However, vehicles that receive the least amount of energy are charged to a higher level.

\subsection{Average performance of charging strategies}
\label{macroscopic_performance}

\begin{figure}[ht!]
\centering
\includegraphics[scale=0.3]{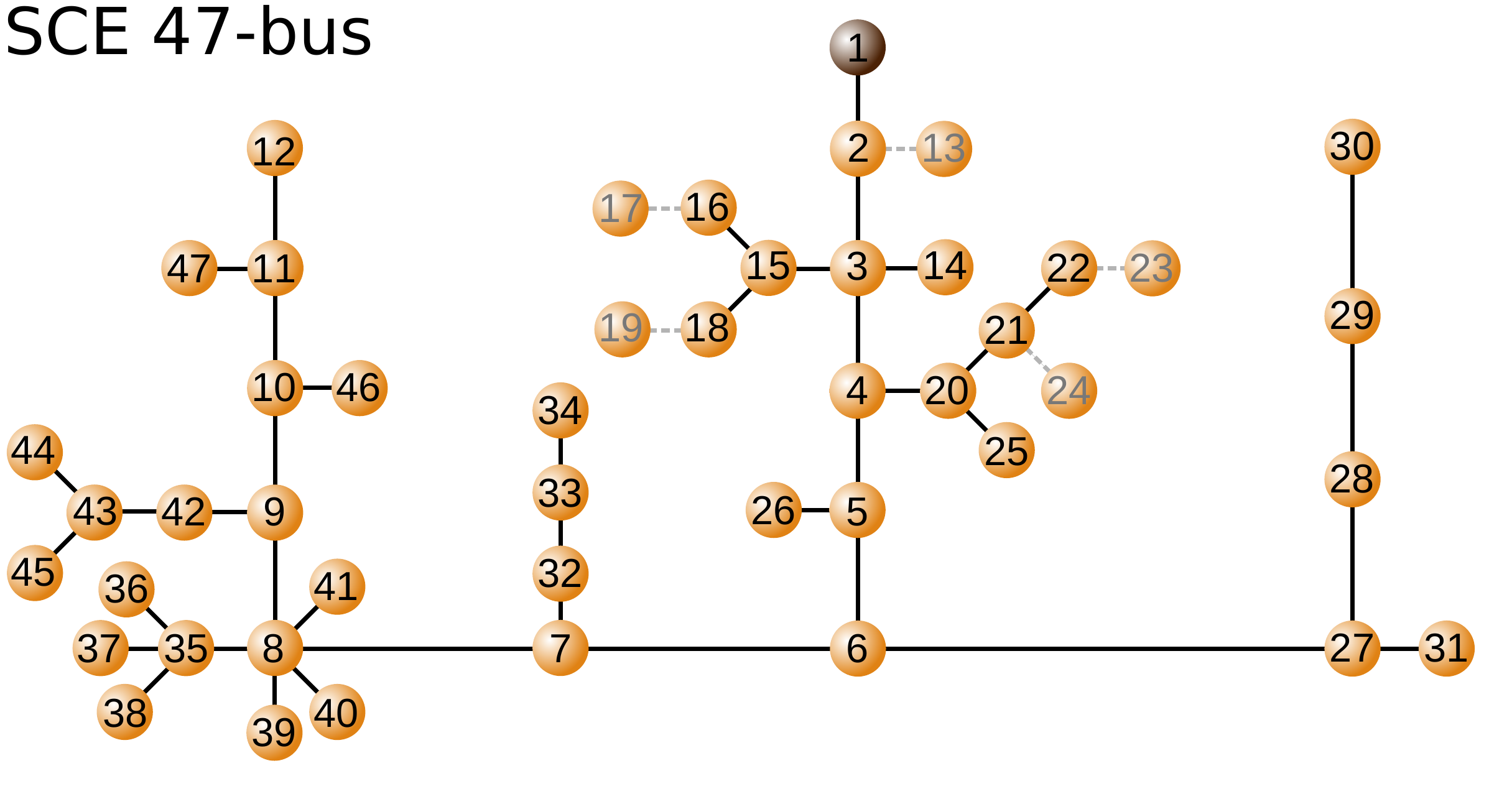}
\caption{Electric network used in numerical experiments exploring the average performance of charging strategies~\cite{Gan_2015}. Node 1 is the root node. Nodes 13, 17, 19, 23 and 24 (in lighter colour) are photovoltaic generators, and we removed them from the network.}
\label{en32}
\end{figure}

To evaluate the average performance of charging strategies, we assume a Poisson arrival process characterized by the arrival rate $\lambda$. The value of $\lambda$ is varied in the range, which is broad enough to observe nearly all vehicles fully charged (free-flow state) and to demonstrate situations when vehicles are charged to only a small degree (congested state). In the numerical experiments, we used the network topology shown in Figure~\ref{en32}, the same setting of parameters as in the previous Section (unless reported otherwise) and we assume a uniform population of strategies. We set the simulation horizon to $10\ 000$ time units and run at least $20$ independent realizations for each value of the arrival rate. 

\begin{figure}[ht!]
\centering
\includegraphics[scale=0.43]{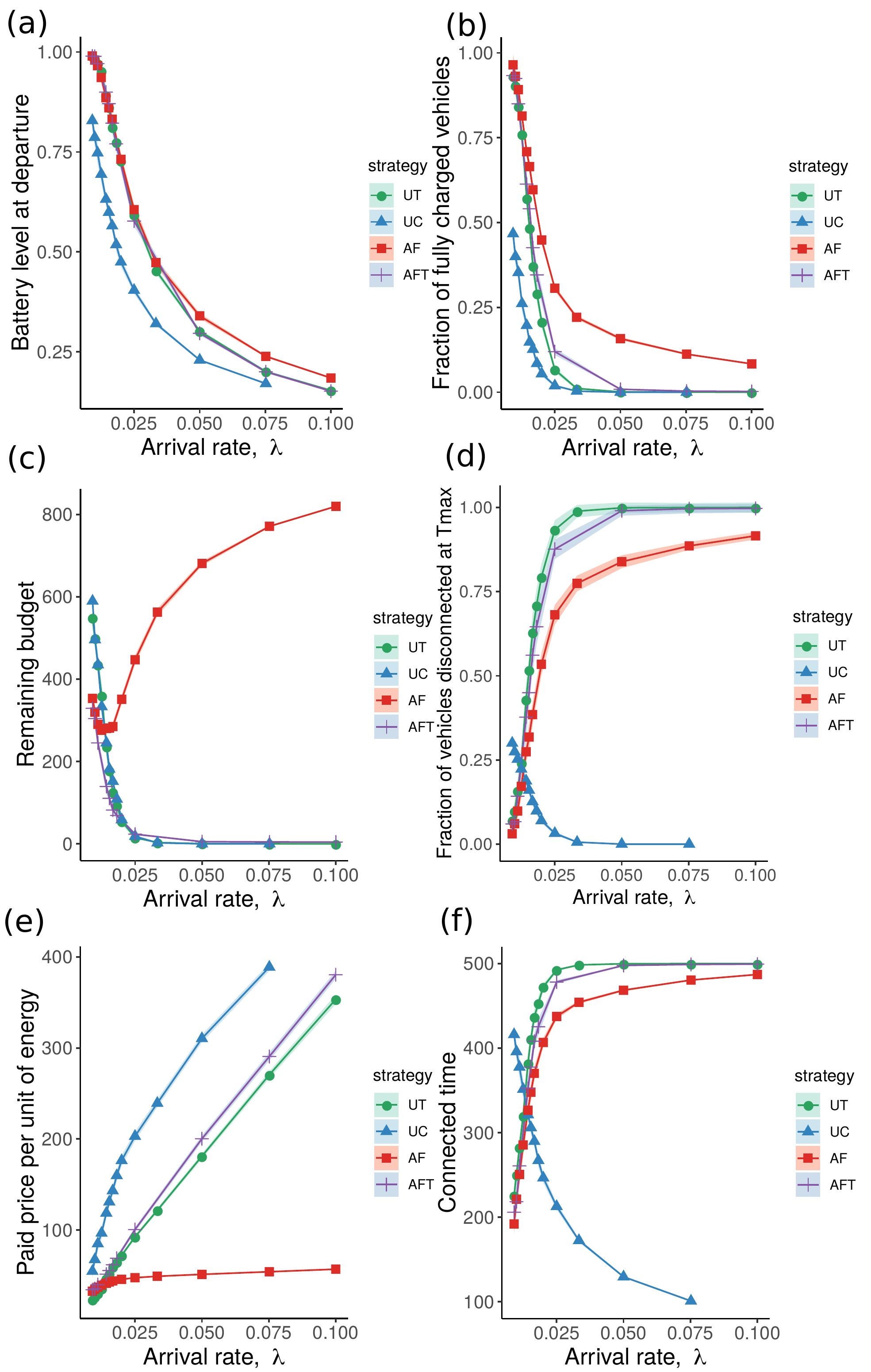}
\caption{Comparison of average values of selected statistics for all investigated charging strategies as a function of arrival rate $\lambda$. Ribbons, which are in most of the cases narrower than the size of symbols, represent 95 percent confidence intervals of average values. Average battery level at departure, number of fully charged vehicles and number of vehicles disconnected at time $T_{max}^l$ have been rescaled.}
\label{fig:performance}
\end{figure}

As expected, the main trend in the Figure~\ref{fig:performance} is that with increasing the arrival rate the state of charge of batteries, measured at the time of departure of the vehicles from the network, decreases. Strategies UT, AF and AFT behave in broadly similar way, while strategy UC reaches significantly smaller average state of charge which is also reflected by the significantly smaller number of vehicles that are fully charged (see Figures~\ref{fig:performance}a-\ref{fig:performance}b). 

In the free-flow state, strategy UC spends the least budget on average (see Figure~\ref{fig:performance}c). However, the vehicles have the smallest battery level when they depart. The average battery level on departure is lower than in other strategies. This is caused by vehicles connected to nodes distant from the root node. Often, they need to increase their willingness to pay parameter to keep the required charging rate and consequently they run out of the budget.

In the congested state, the majority of vehicles following strategies UT, AF and AFT leave the network at the departure time $T_{max}^l$ (see Figure~\ref{fig:performance}d), while vehicles following strategy UC typically run out of the budget. This is caused by the extensive competition in increasing willingness to pay parameter between vehicles following strategy UC and it is also reflected by the high price paid per unit of energy (see Figure~\ref{fig:performance}e). Consequently, the connected time of vehicles following UC strategy is short, when the network is congested and long in free-flow (see Figure~\ref{fig:performance}f). With all other strategies it is other way around. While in the free-flow vehicles are quickly charged and leave the network, and in the congested state they utilize the whole time period they are given for charging and wait for more favourable conditions to collect some energy.

Another interesting observation can be found in Figures~\ref{fig:performance}c and~\ref{fig:performance}e. With increasing arrival rate $\lambda$, the remaining budget typically decreases and the price paid per unit of energy increases, except for strategy AF, where the remaining budget grows and the price paid per unit of energy stays almost constant. This is because in the congested state, the price per unit of received energy is often unaffordably high. In the congested state, significantly more vehicles are fully charged by the AF strategy. This can be attributed to the sharp increase in the willingness to pay parameter when the state of charge  approaches full capacity. The secondary effect of the sharp increase is that other vehicles lower their willingness to pay parameter (to avoid the aggressive user), which increases the likelihood that the vehicle reaches full battery even more.

Strategies AF and  AFT are similar, the main difference is that strategy AFT utilizes information about the departure time of vehicles and strategy AF does not. Thus, it is surprising that strategy AFT charges fully a smaller fraction of vehicles and their average battery level at departure is smaller than for the AF strategy. This can be attributed to indirect effects. When the AFT strategy increases the willingness to pay parameter for a vehicle connected to a more distant node shortly before it leaves the network, this vehicle receives more energy than it would receive with the AF strategy. However, when the network is congested this is at the expense of other vehicles that are closer to the source node. Vehicles positioned closer to the source node consume less network resources and they could receive more energy in total, if using the AF strategy. To summarize, the way strategy AFT uses information about the departure of vehicles $T_{max}^l$, leads to the lower total average charging level of vehicles at the departure, but the vehicles that are connected to more distant nodes from the root node are charged to higher degree. Thus, unexpectedly what the AFT strategy gains compared to the AF strategy in a congested state is not a higher average state of charge of batteries at the time of departure of vehicles but a more even final level of charge across vehicles.

\section{Conclusions}
\label{sec_conclusions}
In this paper, we took some initial steps toward understanding the collective effects potentially emerging in future smart electric grids, arising from approaches that are exploiting flexibility of loads, such as charging of electric vehicles, using basic market mechanisms. The users are allowed to pre-allocate a limited budget, which they wish to spend for a service and they state the time interval within which the service can be provided. We elaborated a mathematical model of the electric distribution network that allocates network resources based on the proportional fairness mechanism. We proposed simple algorithmic strategies without memory capable of handling budget and charging time restrictions using different approaches and we investigated the behaviour emerging from individual strategies in pre-designed situations. Furthermore, we studied the average behaviour of strategies.

Our analyses show that even simple strategies that prescribe how a limited budget is spent can lead to interesting collective effects, such as the emergence of a sorting effect in charging of vehicles, avoidance of aggressive users, or to unexpected indirect effects, when instead of increasing the average state of charge, equality among vehicles is enhanced. Our results vividly illustrate the need to carefully analyse collective effects emerging in complex systems in order to avoid negative effects and to exploit positive ones.

Future steps could focus on the collective behaviour of strategies and identify other interesting behavioural patterns. For example, the number of vehicles could be periodically increasing (congestion builds up) and decreasing (congestion resolves) or scenarios with mixed composition of strategies could be investigated. Interesting results could be obtained by strategies that use learning mechanisms.

\section*{Acknowledgement}
We thank Jonas Knapp for developing early version of charging strategies and numerical experiments. This work was supported by The Alan Turing Institute under the EPSRC grant EP/N510129/1, by the research grants VEGA 1/0463/16 "Economically efficient charging infrastructure deployment for electric vehicles in smart cities and communities",  APVV-15-0179 "Reliability of emergency systems on infrastructure with uncertain functionality of critical elements" and it was facilitated by the FP 7 project ERAdiate [621386] "Enhancing Research and innovation dimensions of the University of \v{Z}ilina in Intelligent Transport Systems".
\bibliographystyle{IEEEtran}
\bibliography{compeng_2018}
\end{document}